# Modelling the Time-dependent VRP through Open Data

Augustin Lombard, Simon Tamayo and Frédéric Fontane

*Abstract*— This paper presents an open data approach to model and solve the vehicle routing problem with time-dependent travel times (TDVRP). The proposed model is based on a multi-layer matrix composed of travel times, replacing the traditional distance matrix. Online cartography services are queried in order to build this matrix. Travel times are obtained for every step in the time discretization. Thus, the model integrates the fact that the travel time between two points is modified during the time horizon. This model is applied to a medium-sized problem in the urban area of Paris using an enhanced Greedy Randomized Adaptive Search Procedure (GRASP). This work intends to build on the current state of the art by proposing a straightforward and open-access method to model and solve the VRP with traffic variability.

## I. INTRODUCTION

THIS paper presents an approach to model the static and deterministic time-dependent vehicle routing problem by using online cartography services information.

The Vehicle Routing Problem (VRP) aims at finding the optimal set of routes to be performed by a fleet of vehicles to serve a given set of customers [1].

The VRP was first introduced in 1959 by Dantzig and Ramser [2] when addressing gasoline delivery optimization problems. The VRP can be solved using exact methods (Fisher et al. [3]; Bard et al. [4]; Toth and Vigo [1]; Chabrier [5]; Qureshi et al. [6]; Azi et al. [7]) or heuristic algorithms (Badeau et al. [8]; Cordeau et Maischberger [9]; Laporte et al. [10]; Jones et al. [11]; Montemanni et al. [12]; Polimeni and Vitettaa [13]; Adulyasak and Jaillet [14]). Nonetheless, the exact approaches have limitations related to the computing time, which makes them less adapted to large problems.

The VRP has many important variants that take into account different delivery configurations and constraints. A quite exhaustive taxonomy of the VRP variants can be found in the work of Eksioglu et al. [15]. One of these variants is the time-dependent VRP also known as TDVRP, which

Manuscript received December 22, 2017; revised January 30, 2018. This work was supported by ADEME (French environment and energy management agency), La Poste, Marie de Paris (City of Paris), Pomona Group and Renault, within the context of the Urban Logistics Chair at Mines ParisTech.

A. Lombard. (phone: +33663013516 e-mail: augustin.lombard@mines-paristech.fr).
S. Tamayo. (e-mail: simon.tamayo@mines-paristech.fr).
F. Fontane. (e-mail: frederic.fontane@mines-paristech.fr).
All authors are with MINES ParisTech, PSL Research University, 60 boulevard Saint-Michel, Paris, France.

allows introducing a notion of "variable cost" in the objective function that is measured in time instead of distance. The TDVRP allows integrating the notion of travel time variability due to road congestion.

The TDVRP was first introduced in 1992 by Malandraki and Daskin [16] and by Hill and Benton [17]. In both cases the problem was formulated as a Mixed Integer Linear Program (MILP) problem and solved through heuristic approaches. The TDVRP has been largely explored in the literature. Gendreau et al. [18] proposed a survey on the problem in which they classify the problem in 3 categories depending on the quality and the evolution of the information related to travel times, namely: 1) Static and deterministic time-dependent VRP; 2) Static and stochastic time-dependent VRP; and 3) Dynamic time-dependent VRP.

One of the limits highlighted in the existing literature is that the time-dependent function (i.e. the evolution of the speed in time) is difficult to obtain, yet it is the basis for a pertinent solution. This paper aims at proposing a model to solve the TDVRP using open-access cartography data to model the time-dependent travel times. More precisely, a static and deterministic time function is represented using a multi-layer travel times matrix, which is built by querying cartography services – such as Google Maps. Thus, with a list of addresses or coordinates, one can find a pseudo-optimum itinerary, taking traffic into account.

This paper is organized as follows. Section 2 states the problem and the research questions. Section 3 presents the proposed approach, the mathematical model and the data collection procedure. Section 4 proposes an application problem in the urban area of Paris and its results. Finally, section 5 concludes with remarks and openings.

## II. PROBLEM STATEMENT

### A. Model the time-dependent VRP with up to date, realistic and open access data

Professional software presently used to create delivery itineraries often avoid taking traffic into account. However, today, traffic is responsible for significant difficulties, such as late deliveries, loss of time for drivers or increased fuel consumption, to name a few. Open access cartography tools such as Google Maps exist to improve a given itinerary and reduce driving times. This paper aims at modeling and optimizing the VRP based on open web services to realistically model traffic during a given period of time.

### B. Research questions

This work aims at formalizing a simple and straight forward model for solving the TDVRP. The proposed

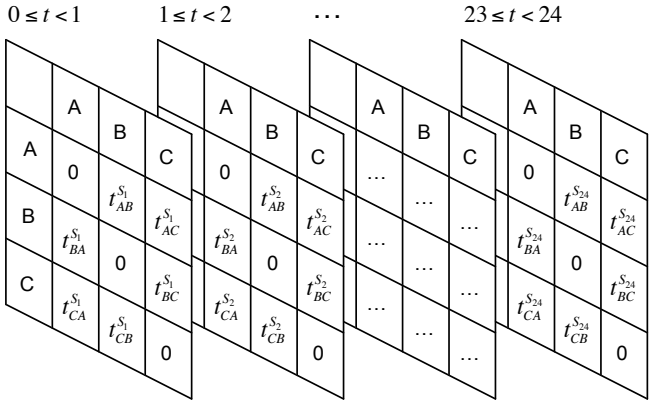

Fig. 1 An example of a multi-layer matrix.

approach should improve the planning of itineraries by taking into account the overall average traffic in a road network. In order to tackle this problem, two research questions must be addressed:

RQ1. How can the time-dependent VRP be modeled and optimized efficiently?

RQ2. How can open-source data be used to model real life time-dependent graphs?

### III. PROPOSED APPROACH

#### A. Multi-layer travel times matrix

In order to model time-dependent travel times, this paper proposes a matrix representing travel times at different moments of a planning horizon; a multi-layer matrix, with one layer for each time-step, with each layer containing the travel times between any points in the problem.

The matrix focuses only in the time needed to go from one point to another (which will depend on traffic). As a result, it is assumed that each layer is not necessarily symmetrical.

Fig. 1 shows an example of the proposed multi-layer matrix.

As shown in Fig. 1 each layer is a matrix containing the average travel times between the nodes of the problem for a time-step. For instance, $t_{AB}^{S_1}$ in the first layer indicates the average travel time from A to B in the first hour of the day.

#### B. Mathematical model

The time-dependent VRP can be modelled as an oriented graph with a node set $V$, and arc set $A$. The first node corresponds to the depot and the rest to the clients. Each arc $(i, j)$ has an associated time-dependent cost $C_{ij}(k_i)$, equal to the travel time between nodes $i$ and $j$ at a departure time $k_i$. In any given time-step, cost respects the triangular inequality. A given route is modelled as a list $R$ of nodes; $R_a$ is the $a^{th}$ node served in the route. $P_i$ is the position of node $i$ in a given route $R$. $T$ represents the length of the planning horizon; and $u_i$ is a dummy variable used to express condition (1d). $x_{ij} \in \{0,1\}$ indicates if in the solution, the vehicle travels from node $i$ to node $j$.

In the model, drop-off and handling times are not considered. Thus, arrival time to a node is also the departure time from the same node. The vehicle must start its tour from the depot at the beginning of the planning horizon ($k_0 = 0$).

With the previous definitions, the problem can be modelled as follows:

$$min \sum_{i \in V} \sum_{\substack{j \in V \\ j \neq i}} x_{ij} C_{ij}(k_i) \tag{1a}$$

$Under$:

$$\forall i \in V, \sum_{j \in V} x_{ij} = 1 \tag{1b}$$

$$\forall j \in V, \sum_{i \in V} x_{ij} = 1 \tag{1c}$$

$$u_i - u_j + n\, x_{ij} \leq n - 1 \quad 2 \leq i, j \leq n, i \neq j \tag{1d}$$

$$\forall i, j \in V, x_{ij} \in \{0,1\} \tag{1e}$$

$$\forall i \in V, u_i \in \mathbb{Z} \tag{1f}$$

$$\forall i \in V, k_i \in [0, T] \tag{1g}$$

Constraints (1b) and (1c) ensure flow conservation: there is exactly one arrival and one departure from every node. Constraint (1d) ensures that every route passes through node 1. Combined with the flow conservation, it proves that all nodes are served in the same route (i.e. there is no subroute). Constraint (1d) comes from the integer formulation of the travelling salesman problem of Miller et al. [19]. It can be proven that $P_i$ are values for the dummy variables $u_i$ that verify constraint (1d).

In a route $R$, node $i$ is at position $P_i$. The nodes visited before node $i$ are nodes $R_a$ with $a = 1, \ldots, P_i-1$. Then, departure time from $i$, $k_i$, is the sum of the time needed to go from node $R_a$ to $R_{a+1}$, with $a = 1, \ldots, P_i-1$. Departure time $k_i$ from a node $i$ can be computed recursively as follows:

$$k_0 = 0 \tag{2a}$$

$$k_i = \sum_{a=1}^{P_i - 1} x_{R_a R_{a+1}} C_{R_a R_{a+1}}(k_{R_a}) \tag{2b}$$

#### C. Open data fetching approach

To create the multi-layer time matrix, we propose using open data. This paper uses Google Maps, and its Distance Matrix API. In order to obtain a time-dependent matrix that represents the average behavior of the modelled zone (i.e. to limit the influence of rare incidents, such as car crashes, etc.), the query must be done in the future. In the present example the query is done in the near future, two weeks after present day.

It is important to note that the time matrix is not symmetrical because of one-way roads, and different levels of congestion.

In terms of distance, the triangular inequality is not necessarily verified because of traffic and the different types of roads (e.g. highway, residential areas). The fastest road is not necessary the shortest. However, the cost in the objective function, the driving time, verifies the triangular inequality. The driving time between two points will always be the smallest.

The matrix is built as follows:
- Select $N$ points on a map representing the clients and the depot.

- Retrieve longitude and latitude coordinates for the previously pinned points.
- Discretize the time horizon in a given number of time steps.
- Query online cartography services and obtain the travel times between all points, for all time steps.
- Assemble the different layers into a single matrix.

*D. Optimization procedure*

This paper uses a greedy randomized adaptive search procedure (GRASP) and an enhancement heuristic to optimize the model. This method builds on the work of Reyes et al. [20], that solved a version of the VRP, the VRP with roaming delivery locations (VRPRDL).

In the first step of the optimization, a GRASP is conducted. $N_{grasp}$ different routes will be created to serve all of the nodes and the best will be kept. Routes are created as follows. For each node not yet inserted on the road, cost of its insertion at any place is calculated. The insertions are sorted in ascending order and the $K_{grasp}$ smallest are kept. One insertion is randomly chosen and the corresponding node is inserted at the corresponding location of the route. When all the nodes are placed, a total cost is calculated. At the end of the GRASP, among the $N_{grasp}$ routes, the one with the lowest cost is kept. The output of the GRASP is a single solution (the best one amongst the $N_{grasp}$ trials).

After the GRASP, an "insertion-deletion" heuristic is used to improve the current solution. $L$ nodes are deleted and then reinserted, following a first-deleted, first-reinserted rule. First, we build a list containing the nodes that would lead to the greatest cost saving in the solution if they were to be deleted. The size of the list is stated by a given parameter $K_{del}$, then the deleted nodes are randomly chosen one at a time and the list is updated after each deletion. Nodes are reinserted in the order they were deleted, at a place chosen randomly in the list of $K_{ins}$ places leading to the lowest cost. This heuristic is repeated $N'$ times.

If parameters $K_{grasp}$, $K_{del}$ and $K_{ins}$ are large numbers, the stochastic search of the algorithm increases, which could help avoid local optima. Indeed, higher stochasticity would need more search iterations to ensure convergence. However, if the algorithm is intended to run during a given computation time, high values for these parameters lead to less robust solutions because of the uncertainty of the insertions.

## IV. APPLICATION

The approach is tested on the urban area of Paris. In the recent years, through-roads along the river Seine were closed for traffic and opened for pedestrians. This resulted in an increase of traffic on the ring road and other roads to cross the city [21]. Managing deliveries with traffic is now a key issue for delivery services.

The example is composed of 30 clients, 1 depot, all in Paris and its urban area. The multi-layer time matrix is queried for the December 12th 2017, from 8 a.m. to 8 p.m. with time steps of 2 hours. As a result, the multi-layer matrix is composed of 16 31x31 matrices. Fig. 2 shows clients and depot on a map of Paris urban area and table I shows the latitude and longitude of all the chosen clients.

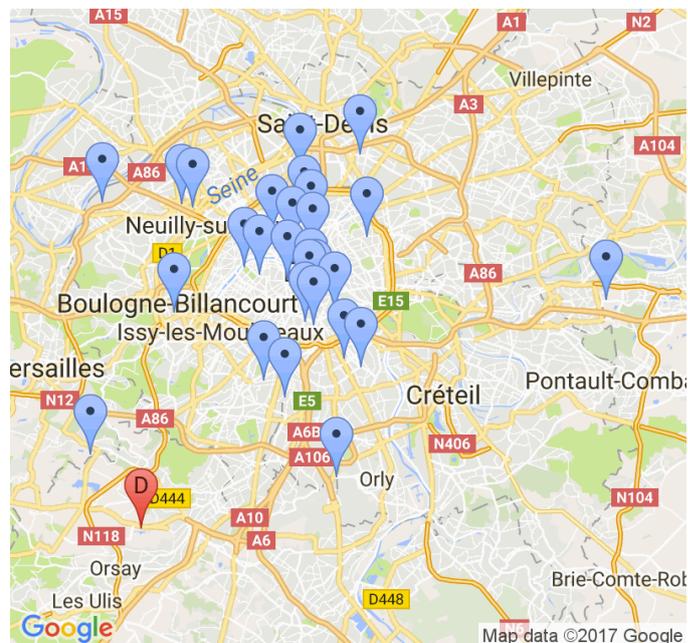

Fig. 2. Map of Paris urban area with 30 clients and one depot (D).

TABLE I
COORDINATES OF THE CLIENTS

| Client index | Latitude (°) | Longitude (°) |
|---|---|---|
| 0 (Depot) | 48.845761 | 2.339546 |
| 1 | 48.847397 | 2.348344 |
| 2 | 48.758064 | 2.169373 |
| 3 | 48.841295 | 2.588015 |
| 4 | 48.841575 | 2.347142 |
| 5 | 48.858862 | 2.294403 |
| 6 | 48.886944 | 2.343126 |
| 7 | 48.745642 | 2.369067 |
| 8 | 48.854819 | 2.306285 |
| 9 | 48.870349 | 2.33312 |
| 10 | 48.892779 | 2.243434 |
| 11 | 48.890971 | 2.252033 |
| 12 | 48.893764 | 2.178884 |
| 13 | 48.834553 | 2.237409 |
| 14 | 48.788931 | 2.327215 |
| 15 | 48.797628 | 2.309662 |
| 16 | 48.809578 | 2.375374 |
| 17 | 48.805202 | 2.388838 |
| 18 | 48.827828 | 2.350659 |
| 19 | 48.82934 | 2.34343 |
| 20 | 48.834854 | 2.367539 |
| 21 | 48.841673 | 2.351875 |
| 22 | 48.852078 | 2.329325 |
| 23 | 48.867007 | 2.349496 |
| 24 | 48.876527 | 2.316532 |
| 25 | 48.879352 | 2.348189 |
| 26 | 48.909403 | 2.339389 |
| 27 | 48.875115 | 2.393429 |
| 28 | 48.919156 | 2.388138 |
| 29 | 48.917793 | 2.243537 |
| 30 | 48.718288 | 2.210583 |

The matrix was generated by querying Google Maps API [22] using the parameters mode = "driving" and traffic_model = "best_guess".

*A. Parameters*

In the first phase of the optimization, the GRASP, $N = 30$ routes were created with each the 30 nodes. The best route was kept. The insertion was chosen amongst the $K = 3$ best insertions.

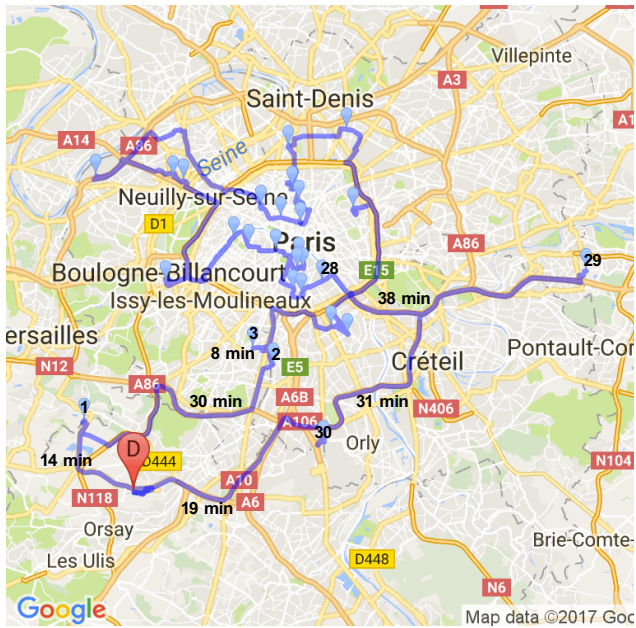
Fig. 3 Optimization results with the classical VRP.

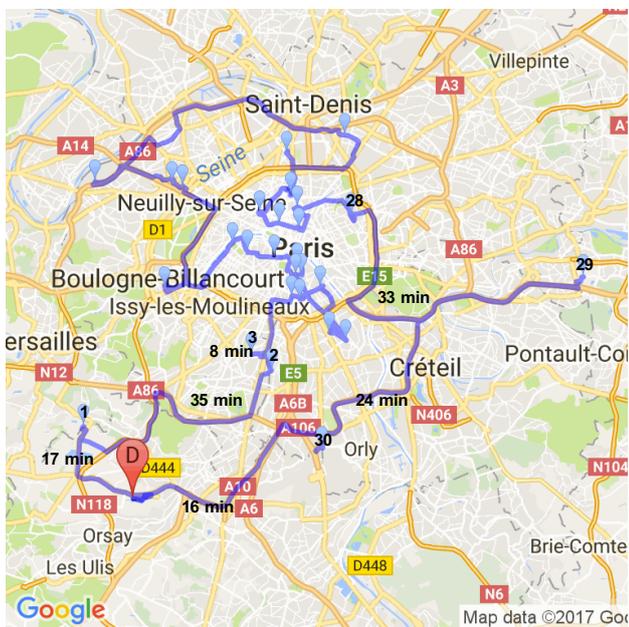
Fig. 4 Optimization results with the TDVRP.

In the second phase of the optimization, the heuristic is run $N' = 20$ times, removing $L = 6$ nodes each time. These parameters were chosen as a trade off between efficiency and computation time. It is important to note that improvement in the choice of these parameters will lead to a better solution.

### B. Results

To assess the quality of the optimization process, the best route was plotted on Google Maps to show the real itinerary taken by the vehicle to go from a node to another. Time-dependent VRP is also compared to the classical VRP, where the travel times matrix is the average of the multi-layer matrix along the time steps. To do so, total cost was compared in both cases.

On Fig. 3, the itinerary of the delivery vehicle can be seen. The vehicle starts with customer 1, does his loop, and finishes with customer 30. The itinerary seems rational, even if there are a few crossings. Crossings in the best itineraries can be explained by sub-optimality of the solution or by traffic and road types (e.g. highway versus city streets). It can be faster to exit a highway after the delivery location is passed by, which will imply crossing roads. The general pattern is the same between the two scenarios, with an 8-shape versus an "almost" 8-shape, and slight modifications in the north of the city. The modifications are explained by the traffic. The order the clients are served can change to minimize the time the vehicle is used.

In Fig. 3, the classical VRP, the pseudo-optimum route lasts 7 hours and 28 minutes. In Fig. 4, the pseudo-optimum route lasts 7 hours and 35 minutes. The relative gap is 1.7%, which seems rather small since the driving time can double during the day on some arcs. In the model, service time is neglected. When a vehicle arrives at a node, it leaves immediately. Since the vehicle starts at 8 a.m. from the depot, it has to go through morning congestion. However, in the solution, the vehicle is back before the afternoon congestion.

Fig. 3 and Fig. 4 indicate some of the travel times between nodes at the beginning and at the end of the route. For the first three nodes, the classical VRP is fastest by 8 minutes, whereas at the end of the itinerary it is slower by 15 minutes. This explains why the TDVRP can lead to faster solutions compared to the classical VRP. That is to say, that although the chosen roads are congested, it is faster to use them in the morning because the return roads are even more congested at the same time. Waiting for the return roads to be less congested leads to time savings.

The solution returned by the TDVRP is not necessarily better in terms of cost, because the cost in the classical VRP is computed with the average travel times matrix. Nonetheless it must be highlighted that the classical VRP will indeed encounter surprises related the congestion, that the TDVRP will be certainly avoiding.

We computed the cost (i.e. travel time) of the itinerary presented in Fig. 4 taken backwards. To go from the depot to node 28 (passing through node 30 and 29), takes 74 minutes at 6 p.m. whereas it takes 109 minutes at 8 a.m.

### C. Sensibility of the model

In order to verify the robustness of the algorithm, a set of 20 instances of optimization was run.

Table II presents the result of these 20 simulations with the set of customers and depot proposed above. $C_{ML}$ indicates the cost of the instance using the multi-layer matrix, $C_{2D}$ indicates the cost of the instance where the travel times matrix is the average along the time steps of the multi-layer matrix. Instances differ by the seed given to the pseudo-random number generator.

Results range from 6h59 to 7h53. This is a difference of just less than an hour. This variability comes from the parameters $K_{grasp}$, $K_{del}$ and $K_{ins}$ of the optimization method, which allows a non-optimal insertion or deletion of a node in a route, to get out of local minima. There is a great sensibility to these parameters. If $K_{grasp}$, $K_{del}$ and $K_{ins}$ were fixed to 1, then the optimization would be fully deterministic, and the results would all be the same. The variability would be reduced by performing more steps of the GRASP or "insertion-deletion" method. If a very large number of solutions was generated in the GRASP method, the solution would be expected to be the same in each column of the table. To be able to have 20 instances tested, each with the multi-layer and the 2D matrixes, in a reasonable time, 30 GRASP and 20 insertions-deletions are done. This is not a search of the global-optimum.

TABLE II
PERFORMANCE ANALYSIS

| Instance | $C_{ML}$ (h :min) | $C_{2D}$ (h :min) | $\frac{C_{ML}-C_{2D}}{C_{2D}}$ % |
|---|---|---|---|
| 1 | 7 :56 | 7 :24 | 7,2 |
| 2 | 7 :43 | 7 :46 | -0.60 |
| 3 | 7 :02 | 7 :50 | -10 |
| 4 | 7 :43 | 7 :45 | -0.35 |
| 5 | 7 :22 | 7 :00 | 5.2 |
| 6 | 7 :21 | 7 :35 | -3.1 |
| 7 | 7 :03 | 7 :52 | -10 |
| 8 | 7 :08 | 7 :36 | -6.1 |
| 9 | 7 :53 | 7 :44 | 2.1 |
| 10 | 7 :05 | 7 :34 | -6.5 |
| 11 | 7 :10 | 7 :38 | -6.2 |
| 12 | 7 :46 | 7 :21 | 5.6 |
| 13 | 7 :31 | 7 :14 | 3.8 |
| 14 | 7 :10 | 7 :35 | -5.4 |
| 15 | 7 :47 | 7 :31 | 3.4 |
| 16 | 7 :13 | 7 :37 | -5.4 |
| 17 | 7 :41 | 7 :12 | 6.7 |
| 18 | 7 :24 | 7 :41 | -3.8 |
| 19 | 7 :31 | 7 :29 | 0.52 |
| 20 | 6 :59 | 7 :30 | -6.9 |
| Mean | 7 :25 | 7 :33 | -4.3 |

Sometimes, the multi-layer matrix leads to better solutions than the 2D matrix. On average, the multi-layer matrix is better by 4.3%. Taking traffic variability into account helps building better solutions.

### D. Technical considerations in the proposed approach

The numerical method and optimizer were programmed using Python. The GRASP was based on the pseudo-code presented by Reyes et al. [20]. No optimizer module was used.

Retrieving data from Google Maps limits the size of the test problem, since there are quotas for the number of queries per day. For free, one can make 2,500 queries per day per project. The limit is 100,000 per project when paying. If traffic information is required, no more than 100 queries at a time are allowed. A multi-layer matrix of 24 time-steps and $N$ nodes will require $24*N^2$ queries. The biggest problem that can be queried in a day has 64 nodes. Because of the 100 queries at a time limit, the method has to be well structured.

Concerning calculation time, the machine has a 2.5 GHz Intel Core i5. One core is used. The 30 GRASP iterations take 80 seconds on average, and the 20 insertions-deletions 10 seconds.

## V. CONCLUSION AND PERSPECTIVES

This paper proposed a new approach to model and optimize the time-dependent vehicle routing problem (TDVRP). A discretization of time was used to create a multi-layer matrix, in which each layer contains the travel times matrix for the nodes of the problem at a given time of the day. One of the originalities of this work yields in the proposed mathematical model, in which the costs of arcs depend on the time at which the vehicle passes through them. As a result, when minimizing overall cost, optimization algorithms will favor solutions in which the "most expensive arcs" are taken in their "less expensive" times. In other words, if a given delivery point is located in a road that is heavily congested between 8 a.m. and 9 a.m., the optimization logic will avoid it until 9 a.m., creating solutions that perform other deliveries before going to the aforementioned point, positively avoiding the congested road at the most congested time.

The proposed model presents a concrete example that can be applied to any city, of planning deliveries taking traffic into account. Yet there are several improvement perspectives, especially to integrate more realistic traffic considerations. Possible improvements include parallelizing the code to compute more difficult problems more efficiently or to add complexity to the problem, such as limited capacity, stochasticity, time windows or roaming delivery locations.

Open data cartography services, such as Google maps, offer accessible data to build realistic models. The proposed approach is based in a discretization of time. It is important to consider that this choice leads to faster computation times but it also implies a rough approach from the traffic-modeling standpoint. An interesting perspective would be to create a multi-layer matrix with a denser discretization, which could be used to interpolate the time horizon and thus, obtain a continuous time-dependency function.

REFERENCES

[1] P. Toth and D. Vigo, *The vehicle routing problem*, vol. 9. 2002.
[2] G. B. Dantzig and J. H. Ramser, "The Truck Dispatching Problem," *Manage. Sci.*, vol. 6, no. 1, pp. 80–91, 1959.
[3] M. L. Fisher, K. O. Jornsten, and O. B. G. Madsen, "Vehicle routing with time windows: Two optimization algorithms," *Oper. Res.*, vol. 45, no. 3, pp. 488–492, 1997.
[4] J. Bard, L. Huang, M. Dror, and P. Jaillet, "A branch and cut algorithm for the VRP with satellite facilities," *IEE Trans.*, no. 30 (9), pp. 821–834, 1998.
[5] A. Chabrier, "Vehicle Routing Problem with elementary shortest path based column generation," *Comput. Oper. Res.*, vol. 33, no. 10, pp. 2972–2990, 2006.
[6] A. G. Qureshi, E. Taniguchi, and T. Yamada, "An exact solution approach for vehicle routing and scheduling problems with soft time windows," *Transp. Res. Part E Logist. Transp. Rev.*, vol. 45, no. 6, pp. 960–977, 2009.
[7] N. Azi, M. Gendreau, and J. Y. Potvin, "An exact algorithm for a single-vehicle routing problem with time windows and multiple routes," *Eur. J. Oper. Res.*, vol. 178, no. 3, pp. 755–766, 2007.
[8] P. Badeau, F. Guertin, M. Gendreau, J.-Y. Potvin, and E. Taillard, "A parallel tabu search heuristic for the vehicle routing problem with time windows," *Transp. Res. Part C Emerg. Technol.*, vol. 5, no. 2, pp. 109–122, 1997.
[9] J.-F. Cordeau and M. Maischberger, "A Parallel Iterated Tabu Search Heuristic for Vehicle Routing Problems," *Comput. Oper. Res.*, vol. 39, pp. 2033–2050, 2011.
[10] G. Laporte, M. Gendreau, J.-Y. Potvin, and F. Semet, "Classical and modern heuristics for the vehicle routing problem," *Int. Trans. Oper. Res.*, vol. 7, no. 4–5, pp. 285–300, 2000.
[11] D. F. Jones, S. K. Mirrazavi, and M. Tamiz, "Multi-objective meta-heuristics: An overview of the current state-of-the-art," *Eur. J. Oper. Res.*, vol. 137, no. 1, pp. 1–9, 2002.
[12] R. Montemanni, L. M. Gambardella, A. E. Rizzoli, and A. V. Donati, "Ant colony system for a dynamic vehicle routing problem," *J. Comb. Optim.*, vol. 10, no. 4, pp. 327–343, 2005.
[13] A. Polimeni and A. Vitetta, "An Approach for Solving Vehicle Routing Problem with Link Cost Variability in the Time," *Procedia - Soc. Behav. Sci.*, vol. 39, pp. 607–621, 2012.

[14] Y. Adulyasak and P. Jaillet, "Models and Algorithms for Stochastic and Robust Vehicle Routing with Deadlines," *Transp. Sci.*, vol. 50, no. 2, pp. 608–626, 2016.

[15] B. Eksioglu, A. V. Vural, and A. Reisman, "The vehicle routing problem: A taxonomic review," *Comput. Ind. Eng.*, vol. 57, no. 4, pp. 1472–1483, 2009.

[16] C. Malandraki and M. S. Daskin, "Time Dependent Vehicle Routing Problems: Formulations, Properties and Heuristic Algorithms," *Transp. Sci.*, vol. 26, no. 3, pp. 185–200, 1992.

[17] A. V Hill and W. C. Benton, "Modelling Intra-City Time-Dependent Travel Speeds for Vehicle Scheduling Problems," *J. Oper. Res. Soc.*, vol. 43, no. 4, pp. 343–351, Apr. 1992.

[18] M. Gendreau, G. Ghiani, and E. Guerriero, "Time-dependent routing problems: A review," *Comput. Oper. Res.*, vol. 64, pp. 189–197, 2015.

[19] C. E. Miller, A. W. Tucker, and R. A. Zemlin, "Integer Programming Formulation of Traveling Salesman Problems," *J. ACM*, vol. 7, no. 4, pp. 326–329, Oct. 1960.

[20] D. Reyes, M. Savelsbergh, and A. Toriello, "Vehicle routing with roaming delivery locations," *Transp. Res. Part C Emerg. Technol.*, vol. 80, pp. 71–91, 2017.

[21] F. Awada, D. Nguyen-Luong, N. Boichon, M. Bouleau, J. Courel, L. Debrincat, N. Pauget, V. Ghersi, A. Kauffmann, F. Mietlicki, and R. René-Bazin, "Fermeture des voies sur berges rive droite à Paris," 2017.

[22] Google, "Google Maps Distance Matrix API." [Online]. Available: https://google-developers.appspot.com/maps/documentation/distance-matrix/start. [Accessed: 20-Dec-2017].